\documentclass[a4paper,12pt]{article}
\usepackage[cp1251]{inputenc}
\usepackage{amssymb,amsmath,amsthm}
\usepackage[russian]{babel}
\usepackage{indentfirst}
\usepackage{hyperref,hypbmsec}
\usepackage[pdftex]{graphicx}
\numberwithin{equation}{section}

\newtheorem{Th}{\hskip\parindent Theorem}[section]
\newtheorem{Le}{\hskip\parindent Lemma}[section]
\newtheorem{Sl}{\hskip\parindent Corollary}[section]
\newtheorem{Zam}{\hskip\parindent Remark}[section]
\newtheorem{Hyp}{\hskip\parindent Conjecture}[section]
\newcommand{\A}{\mathcal{A}}
\newcommand{\E}{\mathfrak{C}}
\newcommand{\R}{\mathfrak{R}}

\newcommand{\D}{\mathfrak{D}}
\newcommand{\N}{\mathbb{N}}

\newcommand{\TT}{\mathcal{T}}

\newcommand{\M}{\mathfrak{M}}
\newcommand{\MM}{\mathfrak{M}^{*}}
\newcommand{\NN}{\mathfrak{N}}
\newcommand{\1}{\mathbf{1}}

\newcommand{\Om}{\widetilde{\Omega}}
\newcommand{\q}{\mathbf{q}}
\newcommand{\KK}{\overline{K}}
\newcommand{\LL}{\hat{l}}

\newcommand{\rr}{\mathbb{R}}
\newcounter{propet}

\renewcommand{\le}{\leqslant}\renewcommand{\ge}{\geqslant}

\begin{document}
\author{ D.\,A.\,Frolenkov\footnote{Research is supported by RFFI (grant № 12-01-31165)} \quad
I.\,D.\,Kan,\footnote{Research is supported by RFFI  (grant № 12-01-00681-а)} }
\title{
%\begin{flushright}
%\small{Светлой памяти\\ профессора Н.М. Коробова\\ посвящается.}
%\end{flushright}
A reinforcement of the Bourgain-Kontorovich's theorem-ІІ.
}
\date{}
\maketitle
\begin{center}
    Abstract
\end{center}

%\begin{abstract}
Zaremba's conjecture (1971) states that every positive integer number $d$ can be represented as a denominator (continuant) of a finite continued fraction  $\frac{b}{d}=[d_1,d_2,\ldots,d_{k}],$ whose partial quotients $d_1,d_2,\ldots,d_{k}$ belong to a finite alphabet  $\A\subseteq\N.$ In this paper it is proved for an alphabet $\A,$ such that the Hausdorff dimension $\delta_{\A}$ of the set of infinite continued fractions whose partial quotients belong to $\A,$ that the set of numbers $d,$  satisfying Zaremba's conjecture with the alphabet $\A,$ has positive proportion in $\N.$ The result improves our previous reinforcement of the corresponding Bourgain-Kontorovich's theorem.
%\end{abstract}

Bibliography: 8 titles.

\textbf{Keywords:\,} continued fraction, continuant, the circle method, exponential sums.\par
%\begin{keywords}
%\textbf{Keywords:} Frobenius numbers, exponential sums .%\end{keywords}

\renewcommand{\proofname}{{\bf Proof}}

\setcounter{Zam}0
%\newpage
%\tableofcontents
%\newpage
%\part{\large{Введение}}
\section{Introduction}
Let $\R_{A}$ be the set of rational numbers whose continued fraction expansion has all partial quotients being bounded by  $A:$
$$
\R_{A}=\left\{\frac{b}{d}=[d_1,d_2,\ldots,d_{k}]\Bigl| 1\le d_{j}\le A \,\mbox{для}\, j=1,\ldots,k\right\}.
$$
Let $\D_{A}$ be the set of denominators of numbers in $\R_{A}$:
$\D_{A}=\left\{d\Bigl| \exists b:\, (b,d)=1, \frac{b}{d}\in\R_{A}\right\}.$
And set
$$
\D_{A}(N)=\left\{d\in\D_{A}\Bigl| d\le N \right\}.
$$
\begin{Hyp}(\textbf{Zaremba's conjecture} ~\cite[p. 76]{Zaremba},\,1971 ).
For sufficiently large $A$  one has
$$\D_{A}=\N.$$
\end{Hyp}
Bourgain and Kontorovich suggested that the problem should be generalized in the following way. Let $\A \in \N$ be any finite alphabet ($|\A|\ge2$) and let $\R_{\A}$ and $\E_{\A}$  be the set of finite and infinite continued fractions
whose partial quotients belong to $\A.$ And let
$\D_{\A}(N)=\left\{d\Bigl| d\le N,\, \exists b: (b,d)=1,\, \frac{b}{d}\in\R_{\A}\right\}$
be the set of denominators bounded by $N.$ Let $\delta_{\A}$ be the Hausdorff dimension of the set $\E_{\A}.$ In the paper ~\cite{FK} we proved the following theorem on the basis of the method devised by Bourgain-Kontorovich~\cite{BK}.
\begin{Th}\label{uslov1}
For any alphabet $\A$ with
\begin{equation}\label{KFcondition1}
\delta_{\A}>1-\frac{5}{\sqrt{369}+23}=0,8815\ldots,
\end{equation}
the following inequality (positive proportion) holds
\begin{equation}\label{BKresult}
\#\D_{\A}(N)\gg N
\end{equation}
\end{Th}
The main result of the paper is the following theorem.
\begin{Th}\label{uslov2}
For any alphabet $\A$ with
\begin{equation}\label{KFcondition2}
\delta_{\A}>\frac{7}{8}=0,875,
\end{equation}
the inequality \eqref{BKresult} holds.
%("положительная пропорция")
%\begin{equation}\label{KFresult}
%\#\D_{\A}(N)\gg N.
%\end{equation}
\end{Th}
The paper is a sequel of our article ~\cite{FK}. So we will heavily refer to statements and constructions in ~\cite{FK}. It should be mentioned that the proof of the Theorem \ref{uslov2} repeats significantly  the proof of the Theorem \ref{uslov1} in ~\cite{FK}. \par
Throughout $\epsilon_0=\epsilon_0(\A)\in(0,\frac{1}{2500}).$ For two functions $f(x), g(x)$ the Vinogradov notation $f(x)\ll g(x)$  means that there exists a constant $C,$ depending on $A,$ such that $|f(x)|\le Cg(x).$ Also a traditional notation $e(x)=\exp(2\pi ix)$ is used. The cardinality of a finite set $S$ is denoted either $|S|$ or $\#S.$ $[\alpha]$ and $\|\alpha\|$ denote the integral part of $\alpha$ and the distance from $\alpha$ to the nearest integer respectively.
\section{Estimates of exponential sums.}
We define the exponential sum $S_N(\theta)$ as follows
\begin{gather}\label{def Sn 0}
S_N(\theta)=\sum_{\gamma\in\Omega_{N} }e(\theta\|\gamma\|),
\end{gather}
where $\Omega_{N}=\Omega_{N}(\A)$ is a proper set of matrices (ensemble) constructed in ~\cite[глава II]{FK}.We use the following norm $\|\gamma\|=\max\{|a|,|b|,|c|,|d|\}$ for the matrix $\gamma=
\begin{pmatrix}
a & b \\
c & d
\end{pmatrix}$  
One can find some more facts about this norm in ~\cite[§5]{FK}. It was obtained in ~\cite[§7]{FK} that to prove the inequality \eqref{BKresult} it is sufficient to obtain the following estimate
\begin{gather}\label{result}
\int_0^1\left|S_N(\theta)\right|^2d\theta\ll\frac{1}{N}|\Omega_{N}|^2.
\end{gather}
It follows from the Dirichlet's theorem that for any $\theta\in[0,1]$ there exist $a,q\in\N\cup\{0\}$ and $\beta\in\rr$ such that
\begin{gather}\label{14-26}
\theta=\frac{a}{q}+\beta,\;(a,q)=1,\; 0\le a\le q\le N^{1/2},\;\beta=\frac{K}{N},\; |K|\le\frac{N^{1/2}}{q},
\end{gather}
with $a=0$ and $a=q$ being possible if only $q=1.$ The purpose of the following reasonings is a slight modification of the results in ~\cite[§12]{FK}.
\begin{Le}(\cite[Lemma 16.1.]{FK})\label{lemma-17-1}
The following inequality holds
\begin{gather}\label{17-1}
\int_0^1\left|S_N(\theta)\right|^2d\theta\le \frac{1}{N}
\mathop{{\sum}^*}_{0\le a\le q\le N^{1/2}}\int\limits_{|K|\le\frac{N^{1/2}}{q}}
\left|S_N(\frac{a}{q}+\frac{K}{N})\right|^2dK,
\end{gather}
where $\mathop{{\sum}^*}$ means that the sum is taken over $a$ and $q$ being coprime for $q\ge1,$ and $a=0,1$ for $q=1.$
\end{Le}
It follows from the statement of Lemma \ref{lemma-17-1} that we need to know how to estimate the following expression
\begin{gather}\label{12-1}
\frac{1}{N}\mathop{{\sum}^*}_{0\le a\le q\le X}\int\limits_{|K|\le Y}
\left|S_N(\frac{a}{q}+\frac{K}{N})\right|^2dK,
\end{gather}
where $Y$ may depend on $q.$ The following reasonings are similar to \cite[Lemma 26 p.145]{Korobov}. Let take a sufficiently large number $T,$ then
\begin{gather}\label{12-2}
\int\limits_{|K|\le Y}\left|S_N(\frac{a}{q}+\frac{K}{N})\right|^2dK\le
\sum_{|l|\le TY}\int\limits_{l/T}^{(l+1)/T}\left|S_N(\frac{a}{q}+\frac{K}{N})\right|^2dK.
\end{gather}
Hence, in any interval $[l/T,(l+1)/T]$ we have
\begin{gather*}
K=\frac{l}{T}+\lambda,\; 0\le\lambda\le\frac{1}{T},\;
\theta=\frac{a}{q}+\frac{K}{N}=\frac{a}{q}+\frac{l}{TN}+\frac{\lambda}{N},
\end{gather*}
\begin{gather}\label{12-3}
\left|S_N(\theta)-S_N(\frac{a}{q}+\frac{l}{TN})\right|\ll\lambda|\Omega_{N}|\Rightarrow
\left|S_N(\theta)\right|^2\ll\left|S_N(\frac{a}{q}+\frac{l}{TN})\right|^2+\lambda^2|\Omega_{N}|^2.
\end{gather}
So
\begin{gather}\label{12-3-1}
\frac{1}{N}\mathop{{\sum}^*}_{0\le a\le q\le X}\int\limits_{|K|\le Y}
\left|S_N(\frac{a}{q}+\frac{K}{N})\right|^2dK\le\frac{1}{N}
\mathop{{\sum}^*}_{0\le a\le q\le X}\sum_{|l|\le TY}\left(\frac{1}{T}\left|S_N(\frac{a}{q}+\frac{l}{TN})\right|^2+
\frac{|\Omega_{N}|^2}{T^3}\right)
\end{gather}
Choosing $T$ sufficiently large we obtain that the investigation of the expression of the form \eqref{12-1} reduce to the investigation of the quantity
\begin{gather}
\frac{1}{TN}
\mathop{{\sum}^*}_{0\le a\le q\le X}\sum_{|l|\le TY}\left|S_N(\frac{a}{q}+\frac{l}{TN})\right|^2.
\end{gather}
Our next purpose is to modify Lemma 12.4. ~\cite{FK}. We formulate the following theorem for convenience of the reader. Let
\begin{gather}\label{15-8-00}
Q_0=\max\left\{\exp\left(\frac{10^5A^4}{\epsilon_0^2}\right),\exp(\epsilon_0^{-5})\right\}.
\end{gather}
%%%%%%%%%%%%%%%%%%%%%%%%%%%%5
%%%%%%%%%%%%%%%%%%%%%%%%%%%%%%%
\begin{Th}\label{theorem13-4}[Theorem 11.5 ~\cite{FK}]
For any $M^{(1)}$ and $M^{(3)}$ such that
\begin{gather}\label{13-42}
Q_0\le M^{(1)},M^{(3)}\le\frac{N}{Q_0},\;M^{(1)}M^{(3)}< N^{1-\epsilon_0},
\end{gather}
the ensemble $\Omega_N$ can be represented in the form
$\Omega_N=\Omega^{(1)}\Omega^{(2)}\Omega^{(3)},$ such that for any matrices $\gamma_1\in\Omega^{(1)},\,\gamma_2\in\Omega^{(2)},\,\gamma_3\in\Omega^{(3)}$ the following inequalities hold
\begin{gather}\label{13-49}
M^{(1)}\ll\|\gamma_1\|\ll(M^{(1)})^{1+2\epsilon_0},\,
M^{(3)}\ll\|\gamma_3\|\ll(M^{(3)})^{1+2\epsilon_0},
\end{gather}
\begin{gather}\label{13-49-1}
\frac{N}{(M^{(1)}M^{(3)})^{1+2\epsilon_0}}\ll\|\gamma_2\|\ll\frac{N}{M^{(1)}M^{(3)}}.
\end{gather}
\end{Th}
%%%%%%%%%%%%%%%%%%%%%%%%%%%%5
%%%%%%%%%%%%%%%%%%%%%%%%%%%%%%%
We denote $\LL=\max\{1,\frac{|l|}{T}\}$ and
%\begin{gather}\label{14-14}
%\sigma_{N,Z}=\sum_{\theta\in Z}\left|S_N(\theta)\right|,
%\end{gather}
\begin{gather}\label{12-4}
P_{Q_1,Q}^{\kappa_1,\kappa}=\left\{
\theta=\frac{a}{q}+\frac{l}{TN}\;\Bigl|\;(a,q)=1,\;  0\le a\le q,\;Q_1\le q\le Q,\; \kappa_1\le l\le\kappa
\right\}.
\end{gather}
For an arbitrary subset $Z\subseteq P_{Q_1,Q}^{\kappa_1,\kappa},$ we denote
\begin{gather}\label{14-42}
\MM(g_{2})=\left\{ (g^{(1)}_{3},g^{(2)}_{3},\theta^{(1)},\theta^{(2)})\in
\Om^{(3)}\times\Om^{(3)}\times Z^2\Bigl|\,
\eqref{14-43}\, \mbox{и}\, \eqref{14-44}\, \mbox{выполнены}
\right\},
\end{gather}
where $\theta^{(i)}=\frac{a^{(i)}}{q^{(i)}}+\frac{l^{(i)}}{TN}, i=1,2$ and
\begin{gather}\label{14-43}
|g_2g^{(1)}_3\frac{l^{(1)}}{TN}-g_2g_3^{(2)}\frac{l^{(2)}}{TN}|_{1,2}\le \frac{1}{M^{(1)}},
\end{gather}
\begin{gather}\label{14-44}
\|g_2g^{(1)}_3\frac{a^{(1)}}{q^{(1)}}-g_2g_3^{(2)}\frac{a^{(2)}}{q^{(2)}}\|_{1,2}=0.
\end{gather}
We recall that the subscripts "1,2" mean that the property holds for both coordinates. The following lemma can de proved in the same manner as Lemma 12.4. in ~\cite{FK}.
%%%%%%%%%%%%%%%%%%%%%%%%%%%%%%%%%%%%%%%%%%%%%%%%%%%%%%%%%%%%%%%%%%%%%%%%%%%%%%%%%%%%%%%%%%%%%%%%%%%%%%%%%%%%%%%%%%%%%%%%%%%%%%%%%%
%%%%%%%%%%%%%%%%%%%%%%%%%%%%%%%%%%%%%%%%%%%%%%%%%%%%%%%%%%%%%%%%%%%%%%%%%%%%%%%%%%%%%%%%%%%%%%%%%%%%%%%%%%%%%%%%%%%%%%%%%%%%%%%%%%
\begin{Le}\label{lemma-14-4}
Let $Z\subseteq P_{Q_1,Q}^{\kappa_1,\kappa}.$ Let $M^{(1)}$ and $M^{(3)}$ satisfy the condition of the \ref{theorem13-4}, and let the inequality
\begin{gather}\label{14-45}
M^{(1)}>150A^2[q^{(1)},q^{(2)}]\max_{\kappa_1\le l\le\kappa}\LL
\end{gather}
holds for any $\theta^{(1)},\theta^{(2)}\in Z.$ Then the following bound holds
\begin{gather}\label{14-46}
\sum_{\theta\in Z}\left|S_N(\theta)\right|\ll (M^{(1)})^{1+2\epsilon_0}\left|\Omega^{(1)}\right|^{1/2}\sum_{g_{2}\in\Om^{(2)}}
\left|\MM(g_{2})\right|^{1/2}.
\end{gather}
\end{Le}
%%%%%%%%%%%%%%%%%%%%%%%%%%%%%%%%%%%%%%%%%%%%%%%%%%%%%%%%%%%%%%%%%%%%%%%%%%%%%%%%%%%%%%%%%%%%%%%%%%%%%%%%%%%%%%%%%%%%%%%%%%%%%%%%%%
%%%%%%%%%%%%%%%%%%%%%%%%%%%%%%%%%%%%%%%%%%%%%%%%%%%%%%%%%%%%%%%%%%%%%%%%%%%%%%%%%%%%%%%%%%%%%%%%%%%%%%%%%%%%%%%%%%%%%%%%%%%%%%%%%%
%%%%%%%%%%%%%%%%%%%%%%%%%%%%%%%%%%%%%%%%%%%%%%%%%%%%%%%%%%%%%%%%%%%%%%%%%%%%%%%%%%%%%%%%%%%%%%%%%%%%%%%%%%%%%%%%%%%%%%%%%%%%%%%%%%
Let state one more lemma of a general nature that was proved in ~\cite[§12]{FK}. A similar statement was used by S.V.\,Konyagin in ~\cite[ 17]{Konyagin}.
\begin{Le}\label{lemma-14-5}
Let $W$ be a finite subset of the interval $[0,1]$ and let $|W|>10.$ Let $f:W\rightarrow \rr_{+}$ be a function such that, for any subset $Z\subseteq W$ the following bound holds
\begin{gather*}
\sum_{\theta\in Z}f(\theta)\le C_1|Z|^{1/2}+C_2,
\end{gather*}
where $C_1,C_2$ are non-negative constants not depending on the set $Z.$ Then the following estimate holds
\begin{gather}\label{14-49}
\sum_{\theta\in W}f^2(\theta)\ll C_1^2\log|W|+C_2\max_{\theta\in W}f(\theta)
\end{gather}
with the absolute constant in Vinogradov symbol.
\end{Le}
%%%%%%%%%%%%%%%%%%%%%%%%%%%%5
%%%%%%%%%%%%%%%%%%%%%%%%%%%%%%%
%%%%%%%%%%%%%%%%%%%%%%%%%%%%%%%%%%%%%%%%%%%%%%%%%%%%%%%%%%%%%%%%%%%%%%%%%%%%%%%%%%%%%%%%%%%%%%%%%%%%%%%%%%%%%%%%%%%%%%%%%%%%%%%%%%
%%%%%%%%%%%%%%%%%%%%%%%%%%%%%%%%%%%%%%%%%%%%%%%%%%%%%%%%%%%%%%%%%%%%%%%%%%%%%%%%%%%%%%%%%%%%%%%%%%%%%%%%%%%%%%%%%%%%%%%%%%%%%%%%%%
%%%%%%%%%%%%%%%%%%%%%%%%%%%%%%%%%%%%%%%%%%%%%%%%%%%%%%%%%%%%%%%%%%%%%%%%%%%%%%%%%%%%%%%%%%%%%%%%%%%%%%%%%%%%%%%%%%%%%%%%%%%%%%%%%%
\section{«The case $\mu=3$.»}
This section corresponds to the section 14 in ~\cite{FK}. So it has the same title. We formulate some results of ~\cite[§14]{FK}
required for proving the estimate \eqref{result}. We also prove a number of lemmas reinforcing the results of ~\cite[§14]{FK}.
The following lemma is a modification of Lemma 14.1. and 14.2. in ~\cite{FK}.  We write $\KK=\max\{1,|K|\}.$
%Следующая лемма является видоизменением леммы 14.1. и леммы 14.2. из ~\cite{FK}. Положим $\KK=\max\{1,|K|\}.$
%\begin{gather}\label{15-8-0}
%Q_0=\max\left\{\exp\left(\frac{10^5A^4}{\epsilon_0^2}\right),\exp(\epsilon_0^{-5})\right\},\quad \KK=\max\{1,|K|\}.
%\end{gather}
%%%%%%%%%%%%%%%%%%%%%%%%%%%%%%%%%%%%
%%%%%%%%%%%%%%%%%%%%%%%%%%%%%%%%%%%%
\begin{Le}\label{lemma-15-1}
If
$\frac{\kappa}{T}Q^{2,5}\le N^{1-\epsilon_0},\, \frac{\kappa_1}{T}\ge Q_0,$
then the following bound holds
\begin{gather}\label{15-1-0}
\frac{1}{T}\sum_{\theta\in P_{Q_1,Q}^{\kappa_1,\kappa}}\left|S_N(\theta)\right|^2\ll
|\Omega_N|^2\left(Q^2\frac{\kappa}{T}\right)^{2\gamma+7\epsilon_0}Q_1^{\gamma-1+5\epsilon_0}.
\end{gather}
\begin{proof}
Let $Z\subseteq P_{Q_1,Q}^{\kappa_1,\kappa}$ be any subset. In the same way as in Lemma 14.1. of ~\cite{FK} we obtain that to satisfy the conditions \eqref{14-43} и \eqref{14-44} it is necessary to have $q^{(1)}=q^{(2)}=\q.$ Then the conditions \eqref{14-43} and \eqref{14-44} can be written as
\begin{gather}\label{15-24}
g_3^{(1)}a^{(1)}\equiv g_3^{(2)}a^{(2)} \pmod{\q},\quad
|g_2(g_3^{(1)}\frac{l^{(1)}}{TN}-g_3^{(2)}\frac{l^{(2)}}{TN})|_{1,2}\le \frac{1}{M^{(1)}}.
\end{gather}
We fix $\theta^{(1)}$ (that is $a^{(1)},\q,l^{(1)}$ are fixed) for which there are $|Z|$ choices. After this we estimate the number of solutions of  \eqref{15-24} independently of $\theta^{(1)}.$ Then (see Lemma 14.8. in ~\cite{FK})  $a^{(2)}$ is uniquely determined and ${x_1y_2\equiv x_2y_1\pmod{\q},}$ where
\begin{gather*}
g_3^{(1)}=(x_1,x_2)^{t},\,g_3^{(2)}=(y_1,y_2)^{t}.
\end{gather*}
In view of the Theorem \ref{theorem13-4}, we obtain from \eqref{15-24} that
\begin{gather*}
\#l^{(2)}\ll\frac{TN}{M^{(1)}\|g_2g_3^{(2)}\|}\ll T(M^{(1)})^{2\epsilon_0}.
\end{gather*}
%Представим $Z$ в следующем виде $Z=\bigcup\limits_{z_1\in Z_1}Z_2(z_1)$, где $z_1=\frac{a}{q},z_2=\frac{l}{TN}, z_1+z_2\in Z.$
Hence,
\begin{gather}\label{15-25}
|\MM(g_2)|\ll|Z|T(M^{(1)})^{2\epsilon_0}
\sum_{g_3^{(1)},g_3^{(2)}\in\Omega^{(3)}}\1_{\{x_1y_2\equiv x_2y_1\pmod{\q}\}}
\end{gather}
We put
\begin{gather}\label{15-10}
M^{(1)}=150A^2Q^2\frac{\kappa}{T},\quad M^{(3)}=\frac{Q_1^{1/2-2\epsilon_0}}{80A^3}.
\end{gather}
Then the congruences in \eqref{15-25} turn into equations and we obtain $|\MM(g_2)|\ll|Z|T(M^{(1)})^{2\epsilon_0}|\Omega^{(3)}|.$ Applying Lemma  \ref{lemma-14-4}, we have
\begin{gather}\label{15-9}
\sum_{\theta\in Z}\left|S_N(\theta)\right|\ll (M^{(1)})^{1+2\epsilon_0}\left|\Omega^{(1)}\right|^{1/2}\sum_{g_{2}\in\Om^{(2)}}|\MM(g_2)|^{1/2}\ll
\left|\Omega^{(1)}\right|^{1/2}\left|\Omega^{(2)}\right|\left|\Omega^{(3)}\right|^{1/2}
|Z|^{1/2}T^{1/2}(M^{(1)})^{1+3\epsilon_0}.
\end{gather}
Hence, using the bound $\left|\Omega^{(i)}\right|\ge(M^{(i)})^{2\delta-\epsilon_0},$ proved in ~\cite[(11.63)]{FK}, we obtain
\begin{gather*}
\sum_{\theta\in Z}\left|S_N(\theta)\right|\ll
|\Omega_N|\frac{(M^{(1)})^{1+3\epsilon_0}}{(M^{(1)}M^{(3)})^{\delta-\epsilon_0/2}}|Z|^{1/2}T^{1/2}.
\end{gather*}
Applying Lemma\ref{lemma-14-5}, we have
\begin{gather*}
\sum_{\theta\in P_{Q_1,Q}^{\kappa_1,\kappa}}\left|S_N(\theta)\right|^2\ll
|\Omega_N|^2T\frac{(M^{(1)})^{2\gamma+7\epsilon_0}(M^{(3)})^{2\gamma+\epsilon_0}}{(M^{(3)})^2}.
\end{gather*}
Using \eqref{15-10}, we obtain \eqref{15-1-0}. Lemma is proved.
\end{proof}
\end{Le}
%%%%%%%%%%%%%%%%%%%%%%%%%%%%%%%%%%%%
%%%%%%%%%%%%%%%%%%%%%%%%%%%%%%%%%%%%
We denote
\begin{gather}\label{14-27}
P_{Q_1,Q}^{(\beta)}=\left\{
\theta=\frac{a}{q}+\beta\;\Bigl|\;(a,q)=1,\;  0\le a\le q,\;Q_1\le q\le Q
\right\}.
\end{gather}
%%%%%%%%%%%%%%%%%%%%%%%%%%%%5
%%%%%%%%%%%%%%%%%%%%%%%%%%%%%%%
%Для каждого $q$ из $Q_1\le q\le Q$ определим каким-либо способом число $a_q,$ такое что $(a_q,q)=1,\,  0\le a_q\le q.$ %Обозначим
%\begin{gather}\label{15-0}
%Z^{*}=\left\{
%\theta=\frac{a_q}{q}+\beta\;\Bigl|\;Q_1\le q\le Q
%\right\}.
%\end{gather}
%Следующая лемма является видоизменением леммы 14.5. из ~\cite{FK}.
%%%%%%%%%%%%%%%%%%%%%%%%%%%%%%%%%%%%
%%%%%%%%%%%%%%%%%%%%%%%%%%%%%%%%%%%%
\begin{Le}\label{lemma-15-4}(see ~\cite[Lemma 14.5.]{FK})
Let the following inequalities hold
\begin{gather*}
N^{\epsilon_0/2}\le Q^{1/2}\le Q_1\le Q,\,\KK Q\le N^{\alpha},
\end{gather*}
and $\alpha\le\frac{1}{2}+\epsilon_0.$ Then for any $Z\subseteq P_{Q_1,Q}^{(\beta)}$ the following bound holds
%Тогда имеет место оценка
%\begin{gather}\label{15-27}
%\sum_{\theta\in P_{Q_1,Q}^{(\beta)}}|S_N(\theta)|\ll
%|\Omega_N|\left(
%N^{1/2+\alpha/2-\delta+3\epsilon_0}Q+
%N^{1-\delta+3\epsilon_0}\frac{Q^{1/2}}{\KK}\right).
%\end{gather}
%Кроме того, для любого $Z\subseteq P_{Q_1,Q}^{(\beta)}$ имеют место оценки
%\begin{gather}\notag
%\sum_{\theta\in Z}|S_N(\theta)|\ll
%|\Omega_N||Z|^{1/2}\left(
%\frac{N^{1-\delta+2\epsilon_0}}{\KK Q^{1/2}_1}+
%N^{1-\frac{3-\alpha}{2}\delta+4,5\epsilon_0}+
%N^{\frac{3+\alpha}{4}-\frac{3-\alpha}{2}\delta+4,5\epsilon_0}Q^{1/2}\right)+\\+
%|\Omega_N|N^{\frac{1+\alpha}{2}-\delta+2,5\epsilon_0}Q.\label{15-27-1}
%\end{gather}
\begin{gather}\label{15-27-2}
\sum_{\theta\in Z}|S_N(\theta)|\ll
|\Omega_N||Z|^{1/2}\left(
\frac{N^{1-\delta+3\epsilon_0}}{(\KK Q_1)^{1/2}}+
\frac{N^{1-\frac{3-\alpha}{2}\delta+4,5\epsilon_0}}{\KK^{1/2}}\right)+|\Omega_N|N^{\frac{1+\alpha}{2}-\delta+2,5\epsilon_0}Q.
\end{gather}
\begin{proof}
%Неравенства \eqref{15-27} и \eqref{15-27-1} были доказаны в ~\cite[лемма 14.5.]{FK}.
The inequality \eqref{15-27-2} is proved in the same manner as  ~\cite[(14.28)]{FK} with the use of Lemma \ref{lemma-15-7}, which will be proved below, instead of Lemma 14.8 in ~\cite{FK}. This completes the proof of the lemma.
\end{proof}
\end{Le}
%%%%%%%%%%%%%%%%%%%%%%%%%%%%%%%%%%%%
%%%%%%%%%%%%%%%%%%%%%%%%%%%%%%%%%%%%
We denote
\begin{gather}\label{14-30}
\NN(g_{2})=\left\{ (g^{(1)}_{3},g^{(2)}_{3},\theta^{(1)},\theta^{(2)})\in
\Om^{(3)}\times\Om^{(3)}\times Z^2\Bigl|\,
\eqref{14-31}\, \mbox{and}\, \eqref{14-32}\, \mbox{hold}
\right\},
\end{gather}
where
\begin{gather}\label{14-31}
\|g_2g^{(1)}_3\frac{a^{(1)}}{q^{(1)}}-g_2g_3^{(2)}\frac{a^{(2)}}{q^{(2)}}\|_{1,2}\le \frac{74A^2\KK}{M^{(1)}},
\end{gather}
\begin{gather}\label{14-32}
|g_2g^{(1)}_3-g_2g_3^{(2)}|_{1,2}\le\min\left\{
\frac{73A^2N}{M^{(1)}},\; \frac{73A^2N}{M^{(1)}\overline{K}}+\frac{N}{\overline{K}}
\left\|g_2g^{(1)}_3\frac{a^{(1)}}{q^{(1)}}-g_2g_3^{(2)}\frac{a^{(2)}}{q^{(2)}}\right\|_{1,2}
\right\}.
\end{gather}
It was proved in ~\cite[лемме 12.3.]{FK} that
\begin{gather}\label{14-34}
\sum_{\theta\in Z}|S_N(\theta)|\ll
(M^{(1)})^{1+2\epsilon_0}\left|\Omega^{(1)}\right|^{1/2}\sum_{g_{2}\in\Omega^{(2)}}
\left|\NN(g_{2})\right|^{1/2}.
\end{gather}
Let $g_3^{(1)}=(x_1,x_2)^{t},\,g_3^{(2)}=(y_1,y_2)^{t},\mathcal{Y}=x_1y_2-y_1x_2.$ We represent the set $\NN(g_{2})$ as the union of the sets $\M_1,\M_2.$ For the first one we have $\mathcal{Y}=0,$ for the second one $\mathcal{Y}\neq0.$ It was proved in ~\cite[лемме 14.7.]{FK} that $|\M_1|\ll_{\epsilon}Q^{2+\epsilon}|\Omega^{(3)}|.$ The following lemma is a modification of Lemma 14.8 in ~\cite{FK}.
%%%%%%%%%%%%%%%%%%%%%%%%%%%%%%%%%%%%
%%%%%%%%%%%%%%%%%%%%%%%%%%%%%%%%%%%%
\begin{Le}\label{lemma-15-7}(см. ~\cite[Lemma 14.8.]{FK})
Under the hypotheses of Lemma \ref{lemma-15-4} one has
\begin{gather}\label{15-55-0}
|\M_2|\ll|Z|\left(\frac{73A^2N}{M^{(1)}}\right)^2\frac{1}{\KK Q_1}
\left(|\Omega^{(3)}|N^{2\epsilon_0}+Q^{1+\epsilon_0}\right).
\end{gather}
\begin{proof}
To simplify we denote $\TT=\frac{73A^2N}{M^{(1)}}.$ It was proved in ~\cite[лемма 14.8.]{FK} that
\begin{gather}\label{15-59}
|\M_2|\le|Z|\sum_{g_3^{(1)}\in\Omega^{(3)}}\sum_{g_3^{(2)}\in\Omega^{(3)}\atop |g_3^{(1)}-g_3^{(2)}|_{1,2}\le\frac{\TT}{\KK}}\1_{\{x_1y_2\equiv x_2y_1\pmod{\q}\}}.
\end{gather}
Changing the variables $z_1=x_1-y_1,\, z_2=x_2-y_2,$ we obtain
\begin{gather}\label{15-60}
|\M_2|\le|Z|\sum_{g_3^{(1)}\in\Omega^{(3)}}\sum_{|z_{1,2}|\le\frac{\TT}{\KK}}\1_{\{x_1z_2\equiv x_2z_1\pmod{\q}\}}.
\end{gather}
We consider three cases.
\begin{enumerate}
  \item Let $z_1>0,\,z_2>0.$ We fix the vector $g_3^{(1)}\in\Omega^{(3)},$ then $x_1z_2-x_2z_1=j\q.$ Let estimate the amount of $j.$ We have
$$x_1-x_2\frac{\TT}{\KK}\le j\q\le x_1\frac{\TT}{\KK}-x_2$$
and, hence, $\#j\ll\frac{T^2}{\q\KK}+1.$ For a fixed $j$ the solution of the congruence is given by the formulae
\begin{gather*}
z_1=z_{1,0}+nx_1,\;z_2=z_{2,0}+nx_2.
\end{gather*}
In view of $x_2\gg\frac{\TT}{(M^{(1)})^{2\epsilon_0}},$ we have $\#n\ll\frac{(M^{(1)})^{2\epsilon_0}}{\KK}+1.$ Thus,
\begin{gather}\label{15-61}
\sum_{g_3^{(1)}\in\Omega^{(3)}}\sum_{0<z_{1,2}\le\frac{\TT}{\KK}}\1_{\{x_1z_2\equiv x_2z_1\pmod{\q}\}}\ll
|\Omega^{(3)}|\left(\frac{\TT^2}{\q\KK}+1\right)\left(\frac{(M^{(1)})^{2\epsilon_0}}{\KK}+1\right).
\end{gather}
It follows from he conditions of Lemma \ref{lemma-15-4} that $\TT^2>\q\KK,$ so one has
\begin{gather}\label{15-62}
\sum_{g_3^{(1)}\in\Omega^{(3)}}\sum_{0<z_{1,2}\le\frac{\TT}{\KK}}\1_{\{x_1z_2\equiv x_2z_1\pmod{\q}\}}\ll
|\Omega^{(3)}|\frac{\TT^2}{\q\KK}N^{2\epsilon_0}.
\end{gather}
  \item Let $z_1>0,\,z_2<0.$ In the same way as in the previous case we obtain
\begin{gather}\label{15-63}
\sum_{g_3^{(1)}\in\Omega^{(3)}}\sum_{0<-z_2,z_{1}\le\frac{\TT}{\KK}}\1_{\{x_1z_2\equiv x_2z_1\pmod{\q}\}}\ll
|\Omega^{(3)}|\frac{\TT^2}{\q\KK}N^{2\epsilon_0}.
\end{gather}
  \item Let$z_1=0.$ One has
\begin{gather}\label{15-64}
\sum_{g_3^{(1)}\in\Omega^{(3)}}\sum_{|z_2|\le\frac{\TT}{\KK}}\1_{\{x_1z_2\equiv0\pmod{\q}\}}\le
\sum_{g_3^{(1)}\in\Omega^{(3)}}\left(\frac{\TT}{q\KK}(x_1,q)+1\right)\le\\\le
|\Omega^{(3)}|+\frac{\TT^2}{q\KK}\sum_{x_1\le \TT}(x_1,q).
\end{gather}
Next
\begin{gather*}
\sum_{x_1\le \TT}(x_1,q)\le\sum_{d|q}d\left(\frac{\TT}{d}+1\right)\ll q^{\epsilon}\TT+q^{1+\epsilon}
\end{gather*}
and so
\begin{gather}\label{15-65}
\sum_{g_3^{(1)}\in\Omega^{(3)}}\sum_{|z_2|\le\frac{\TT}{\KK}}\1_{\{x_1z_2\equiv0\pmod{\q}\}}\ll
|\Omega^{(3)}|+\frac{\TT^2}{q\KK}\left(q^{\epsilon}\TT+q^{1+\epsilon}\right).
\end{gather}
\end{enumerate}
Using \eqref{15-62},\,\eqref{15-63} and \eqref{15-65} we obtain
\begin{gather}\label{15-66}
|\M_2|\ll|Z|\frac{\TT^2}{q\KK}\left(q^{\epsilon}\TT+q^{1+\epsilon}+|\Omega^{(3)}|N^{2\epsilon_0}\right)\le
|Z|\frac{\TT^2}{q\KK}\left(q^{1+\epsilon}+|\Omega^{(3)}|N^{2\epsilon_0}\right).
\end{gather}
Lemma is proved.
\end{proof}
\end{Le}
%%%%%%%%%%%%%%%%%%%%%%%%%%%%%%%%%%%%
%%%%%%%%%%%%%%%%%%%%%%%%%%%%%%%%%%%%
%%%%%%%%%%%%%%%%%%%%%%%%%%%%%%%%%%%%
%%%%%%%%%%%%%%%%%%%%%%%%%%%%%%%%%%%%
%%%%%%%%%%%%%%%%%%%%%%%%%%%%%%%%%%%%
%%%%%%%%%%%%%%%%%%%%%%%%%%%%%%%%%%%%
Based on Lemma \ref{lemma-15-4} and using Lemma \ref{lemma-14-5} we obtain in the same way as ~\cite[следствию 14.1.]{FK} the following statement
%%%%%%%%%%%%%%%%%%%%%%%%%%%%%%%%%%%%%%%%%%%%%%%%%%%%%%%%%%%%%%%%%%%%%%%%%%%%%%%%%%%%%%%%%%%%%%%%%%%%%%%%%%%%%%%%%%%%%%%%%%%%%%%%%%
%%%%%%%%%%%%%%%%%%%%%%%%%%%%%%%%%%%%%%%%%%%%%%%%%%%%%%%%%%%%%%%%%%%%%%%%%%%%%%%%%%%%%%%%%%%%%%%%%%%%%%%%%%%%%%%%%%%%%%%%%%%%%%%%%%
\begin{Sl}\label{sled15-1}
Under the hypotheses of Lemma \ref{lemma-15-4} one has
%\begin{gather}\label{15-74}
%\mathop{{\sum}^*}_{Q_1\le q\le Q\atop 1\le a\le q }\left|S_N(\frac{a}{q}+\frac{K}{N})\right|^2\ll
%|\Omega_N|^2\left(
%\frac{N^{2(1-\delta)+\frac{\alpha-1}{2}+5\epsilon_0}Q}{(\KK Q_1)^{1-2\epsilon_0}}+
%\frac{N^{2(1-\delta)+4\epsilon_0}Q^{1/2}}{\KK^{2-2\epsilon_0}Q_1^{1-2\epsilon_0}}
%\right),
%\end{gather}
\begin{gather}\label{15-74-0}
\sum_{\theta\in P_{Q_1,Q}^{(\beta)}}|S_N(\theta)|^2\ll
|\Omega_N|^2\left(C_1^2Q^{\epsilon_0}+C'_2\right),
\end{gather}
where
%\begin{gather}\label{15-75}
%C_1=\frac{N^{1-\delta+2\epsilon_0}}{\KK Q^{1/2}_1}+
%N^{1-\frac{3-\alpha}{2}\delta+4,5\epsilon_0}+
%N^{\frac{3+\alpha}{4}-\frac{3-\alpha}{2}\delta+4,5\epsilon_0}Q^{1/2},\quad
%C_2=\frac{N^{\frac{3+\alpha}{2}-2\delta+3,5\epsilon_0}Q}{(\KK Q_1)^{1-2\epsilon_0}}.
%\end{gather}
%или
\begin{gather}\label{15-75-1}
C_1=\frac{N^{1-\delta+3\epsilon_0}}{(\KK Q_1)^{1/2}}+
\frac{N^{1-\frac{3-\alpha}{2}\delta+4,5\epsilon_0}}{\KK^{1/2}},\quad
C'_2=\frac{N^{\frac{3+\alpha}{2}-2\delta+3,5\epsilon_0}Q}{(\KK Q_1)^{1-2\epsilon_0}}.
\end{gather}
%\begin{proof}
%Оценки \eqref{15-74} и \eqref{15-74-0} с параметрами из \eqref{15-75} были доказаны в ~\cite[следствия 14.1 и 14.2.]{FK}. %Параметры \eqref{15-75-1} получаются, если вместо \eqref{15-27-1} использовать \eqref{15-27-2}.  Следствие доказано.
%\end{proof}
\end{Sl}
%%%%%%%%%%%%%%%%%%%%%%%%%%%%%%%%%%%%
%%%%%%%%%%%%%%%%%%%%%%%%%%%%%%%%%%%%
%%%%%%%%%%%%%%%%%%%%%%%%%%%%%%%%%%%%%%%%%%%%%%%%%%%%%%%%%%%%%%%%%%%%%%%%%%%%%%%%%%%%%%%%%%%%%%%%%%%%%%%%%%%%%%%%%%%%%%%%%%%%%%%%%%

%%%%%%%%%%%%%%%%%%%%%%%%%%%%%%%%%%%%%%%%%%%%%%%%%%%%%%%%%%%%%%%%%%%%%%%%%%%%%%%%%%%%%%%%%%%%%%%%%%%%%%%%%%%%%%%%%%%%%%%%%%%%%%%%%%

%%%%%%%%%%%%%%%%%%%%%%%%%%%%%%%%%%%%%%%%%%%%%%%%%%%%%%%%%%%%%%%%%%%%%%%%%%%%%%%%%%%%%%%%%%%%%%%%%%%%%%%%%%%%%%%%%%%%%%%%%%%%%%%%%%
\section{Estimates for integrals of $|S_N(\theta)|^2$.}
The proof of Theorem  \ref{uslov2} is similar to the proof of Theorem \ref{uslov1} in ~\cite{FK}. We will need a number of Lemmas from ~\cite[§16]{FK} which will be presented without proof.
%%%%%%%%%%%%%%%%%%%%%%%%%%%%5
%Напомним, что
%\begin{gather*}
%Q_0=\max\left\{\exp\left(\frac{10^5A^4}{\epsilon_0^2}\right),\exp(\epsilon_0^{-5})\right\}.
%\end{gather*}
%%%%%%%%%%%%%%%%%%%%%%%%%%%%%%%
\begin{Le}\label{lemma-17-2}(~\cite[Lemma 16.2.]{FK})
The following inequality holds
\begin{gather}\notag
\int_0^1\left|S_N(\theta)\right|^2d\theta\le
2Q_0^2\frac{|\Omega_N|^2}{N}+
\frac{1}{N}\mathop{{\sum}^*}_{0\le a\le q\le N^{1/2}\atop q>Q_0}\int\limits_{\frac{Q_0}{q}\le|K|\le\frac{N^{1/2}}{q}}
\left|S_N(\frac{a}{q}+\frac{K}{N})\right|^2dK+\\
\frac{1}{N}\mathop{{\sum}^*}_{0\le a\le q\le Q_0}\int\limits_{\frac{Q_0}{q}\le|K|\le\frac{N^{1/2}}{q}}
\left|S_N(\frac{a}{q}+\frac{K}{N})\right|^2dK
+\frac{1}{N}\mathop{{\sum}^*}_{1\le a\le q\le N^{1/2}\atop q>Q_0}\int\limits_{|K|\le\frac{Q_0}{q}}
\left|S_N(\frac{a}{q}+\frac{K}{N})\right|^2dK\label{17-4}
\end{gather}
\end{Le}
%%%%%%%%%%%%%%%%%%%%%%%%%%%%5
%%%%%%%%%%%%%%%%%%%%%%%%%%%%%%%
It is convenient to use the following notation
\begin{gather}\label{17-7}
\gamma=1-\delta,\quad\xi_1=N^{2\gamma+7\epsilon_0}.
\end{gather}
The second and the third integral in the right side of \eqref{17-4} are estimated in the following lemma.
\begin{Le}\label{lemma-17-3-0}(~\cite[Lemmas 16.3, 16.4. and 16.5]{FK})
For $\gamma<\frac{5}{36}-6\epsilon_0$ and $\epsilon_0\in(0,\frac{1}{2500})$ the following inequalities hold
\begin{gather}\label{17-7-1}
\frac{1}{N}\mathop{{\sum}^*}_{1\le a\le q\le N^{1/2}}\int\limits_{|K|\le\frac{Q_0}{q}}
\left|S_N(\frac{a}{q}+\frac{K}{N})\right|^2dK\ll\frac{|\Omega_N|^2}{N}.
\end{gather}
\begin{gather}\label{17-7-8}
\frac{1}{N}\mathop{{\sum}^*}_{0\le a\le q\le Q_0}\int\limits_{\frac{Q_0}{q}\le|K|\le\frac{N^{1/2}}{q}}
\left|S_N(\frac{a}{q}+\frac{K}{N})\right|^2dK\ll\frac{|\Omega_N|^2}{N}.
\end{gather}
\end{Le}
%%%%%%%%%%%%%%%%%%%%%%%%%%%%5
%%%%%%%%%%%%%%%%%%%%%%%%%%%%%%%
It remains to estimate the first integral in the right side of \eqref{17-4}, that is,
\begin{gather}\label{18-7-10}
\frac{1}{N}\mathop{{\sum}^*}_{0\le a\le q\le N^{1/2}\atop q>Q_0}\int\limits_{\frac{Q_0}{q}\le|K|\le\frac{N^{1/2}}{q}}
\left|S_N(\frac{a}{q}+\frac{K}{N})\right|^2dK
\end{gather}
The following lemmas will be devoted to this. We partition the range of summation and integration over $q,\,K$ into six subareas:
%%%%%%%%%%%%%%%%%%%%%%%%%%%%%%%%%%%%%%%%%%%%%%%%%%%%%%%%%%%%%%%%%%%%%%%%%%%%%%%5
\begin{center}
 %  Requires \usepackage{graphicx}
  \includegraphics[width=450pt,height=350pt]{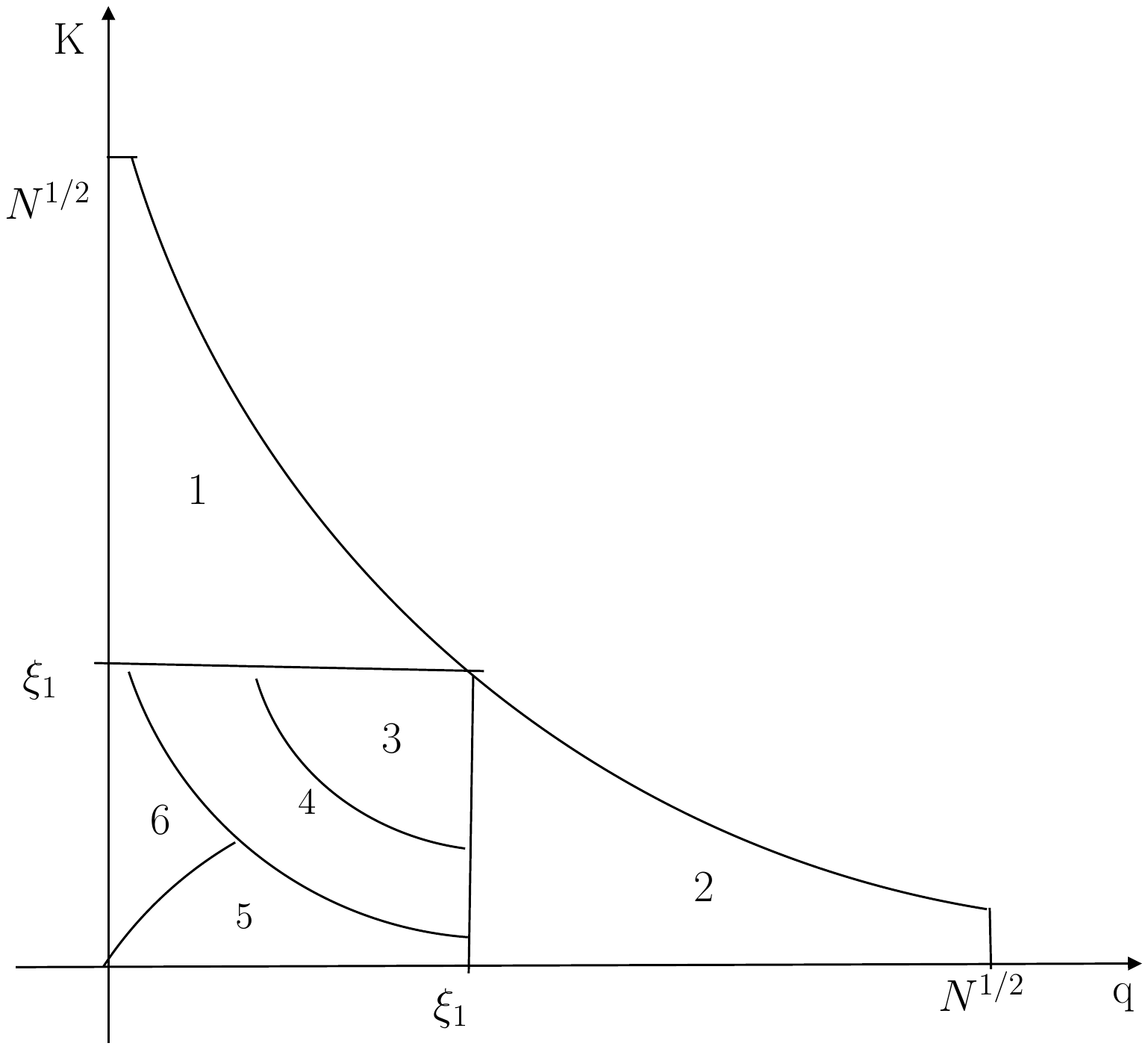}\\
 % \caption{Разбиение на области}\label{ill-1}
\end{center}
Lemma \ref{lemma-17-3} corresponds to the domain 1, Lemma \ref{lemma-17-7} corresponds to the domain 2, Lemma \ref{lemma-17-6-1} corresponds to the domain 3, Lemma \ref{lemma-17-9} corresponds to the domain 4, Lemma \ref{lemma-17-12} corresponds to the domain 5, Lemma \ref{lemma-17-11} corresponds to the domain 6.\par
%%%%%%%%%%%%%%%%%%%%%%%%%%%%%%%%%%%%%%%%%%%%%%%%%%%%%%%%%%%%%%%%%%%%%%%%%%%%%%%%
To prove some lemmas we need the following parameters. Let $N\ge N_{min}=N_{min}(\epsilon_0,\A),$ we denote
\begin{gather}\label{10-2}
J=J(N)=\left[\frac{\log\log N-4\log(10A)+2\log\epsilon_0}{-\log(1-\epsilon_0)}\right],
\end{gather}
where, as usual, $A\ge|\A|\ge2,$ and require the following inequality $J(N_{min})\ge10$ to hold. Now let define a finite sequence $\left\{N_j\right\},$ having set $N_{J+1}=N$ and
\begin{gather}\label{10-8}
N_j=
\left\{
              \begin{array}{ll}
                N^{\frac{1}{2-\epsilon_0}(1-\epsilon_0)^{1-j}}, & \hbox{if $-1-J\le j\le1$;} \\
                N^{1-\frac{1}{2-\epsilon_0}(1-\epsilon_0)^{j}}, & \hbox{if $0\le j\le J$.}
              \end{array}
\right.
\end{gather}
It is obvious that the sequence is well-defined for $j=0$ and $j=1.$ A detailed description of properties of the sequence is given in ~\cite[§9]{FK}.\par
%Леммы \ref{lemma-17-3}~--\ref{lemma-17-5}, следующие далее, также были доказаны в ~\cite[§16]{FK}.

%%%%%%%%%%%%%%%%%%%%%%%%%%%%5
%%%%%%%%%%%%%%%%%%%%%%%%%%%%%%%
\begin{Le}\label{lemma-17-3}(~\cite[Lemma 16.6.]{FK})
%При $\gamma\le\frac{27-\sqrt{633}}{16}-5\epsilon_0,\,\epsilon_0\in(0,\frac{1}{1000})$ имеет место неравенство
The following inequality holds
\begin{gather}\label{17-7-11}
\frac{1}{N}\mathop{{\sum}^*}_{1\le a\le q\le N^{1/2}\atop q>Q_0}\int\limits_{\xi_1\le|K|\le\frac{N^{1/2}}{q}}
\left|S_N(\frac{a}{q}+\frac{K}{N})\right|^2dK\ll\frac{|\Omega_N|^2}{N}.
\end{gather}
\end{Le}
%%%%%%%%%%%%%%%%%%%%%%%%%%%%5
%%%%%%%%%%%%%%%%%%%%%%%%%%%%%%%
%%%%%%%%%%%%%%%%%%%%%%%%%%%%5
%%%%%%%%%%%%%%%%%%%%%%%%%%%%%%%
Let
\begin{gather*}
c_1=c_1(N),\,c_2=c_2(N),\,Q_0\le c_1<c_2\le  N^{1/2},
\end{gather*}
and let
\begin{gather*}
f_1=f_1(N,q),\,f_2=f_2(N,q),\,\frac{Q_0}{q}\le f_1<f_2\le  \frac{N^{1/2}}{q},\\
m_1=\min\{f_1(N,N_j),f_1(N,N_{j+1})\},\,m_2=\max\{f_2(N,N_j),f_2(N,N_{j+1})\}.
\end{gather*}
%%%%%%%%%%%%%%%%%%%%%%%%%%%%5
%%%%%%%%%%%%%%%%%%%%%%%%%%%%%%%
\begin{Le}\label{lemma-17-5}(~\cite[Lemma 16.7.]{FK})
If the functions $f_1(N,q),f_2(N,q)$ are monotonic for $q,$ then the following inequality holds
\begin{gather}\notag
\mathop{{\sum}^*}_{c_1\le q\le c_2\atop 1\le a\le q }\,\int\limits_{f_1\le|K|\le f_2}
\left|S_N(\frac{a}{q}+\frac{K}{N})\right|^2dK\le\\\le
\sum_{j:\,c_1^{1-\epsilon_0}\le N_j\le c_2\,}\int\limits_{m_1\le|K|\le m_2}
\mathop{{\sum}^*}_{N_j\le q\le N_{j+1}\atop 1\le a\le q }\left|S_N(\frac{a}{q}+\frac{K}{N})\right|^2dK.\label{17-7-15}
\end{gather}
\end{Le}
%%%%%%%%%%%%%%%%%%%%%%%%%%%%5
%%%%%%%%%%%%%%%%%%%%%%%%%%%%%%%
\begin{Le}\label{lemma-17-7}
For $\gamma\le\frac{1}{8}-4\epsilon_0,\,\epsilon_0\in(0,\frac{1}{2500})$  the following inequality holds
\begin{gather}\label{17-12}
\frac{1}{N}\mathop{{\sum}^*}_{1\le a\le q\le N^{1/2}\atop q>N^{2\gamma+7\epsilon_0}}
\int\limits_{\frac{Q_0}{q}\le|K|\le\frac{N^{1/2}}{q}}
\left|S_N(\frac{a}{q}+\frac{K}{N})\right|^2dK\ll\frac{|\Omega_N|^2}{N}.
\end{gather}
\begin{proof}
It is sufficient to use Lemma \ref{lemma-17-5} and the estimate \eqref{15-74-0} with \eqref{15-75-1}.
\end{proof}
\end{Le}
%%%%%%%%%%%%%%%%%%%%%%%%%%%%5
%Лемма \ref{lemma-17-6-1} была доказана в ~\cite[§16 лемма 16.10.]{FK}.
%%%%%%%%%%%%%%%%%%%%%%%%%%%%5
\begin{Le}\label{lemma-17-6-1}(~\cite[Lemma 16.10.]{FK})
For $\gamma\le\frac{1}{8}-5\epsilon_0$ the following inequality holds
\begin{gather}\label{17-8-1}
\frac{1}{N}\mathop{{\sum}^*}_{1\le a\le q\le N^{2\gamma+7\epsilon_0}\atop q>N^{\gamma+5\epsilon_0}}
\int\limits_{\frac{N^{3\gamma+12\epsilon_0}}{q}\le|K|\le \xi_1}
\left|S_N(\frac{a}{q}+\frac{K}{N})\right|^2dK\ll\frac{|\Omega_N|^2}{N}.
\end{gather}
\end{Le}
%Лемма \ref{lemma-17-9} была доказана в ~\cite[§16 лемма 16.13.]{FK}.

%%%%%%%%%%%%%%%%%%%%%%%%%%%%5
\begin{Le}\label{lemma-17-9}(~\cite[Lemma 16.13.]{FK})
For $\gamma\le\frac{1}{8}-5\epsilon_0$ the following inequality holds
\begin{gather}\label{17-17}
\frac{1}{N}\mathop{{\sum}^*}_{1\le a\le q\le \xi_1\atop q>Q_0}
\int\limits_{\frac{\xi_1}{q}\le|K|\le\min\{\xi_1,\frac{N^{3\gamma+12\epsilon_0}}{q}\}}
\left|S_N(\frac{a}{q}+\frac{K}{N})\right|^2dK\ll\frac{|\Omega_N|^2}{N}.
\end{gather}
\end{Le}
%%%%%%%%%%%%%%%%%%%%%%%%%%%%5
%Лемма \ref{lemma-17-11}  была доказана в ~\cite[§16 лемма 16.15.]{FK}. Положим $\nu=\frac{3}{2}.$
Let $\nu$ be a positive real number, for example let $\nu\in[1,2].$
%%%%%%%%%%%%%%%%%%%%%%%%%%%%%%%
\begin{Le}\label{lemma-17-11}(~\cite[lemma 16.15.]{FK})
For $\gamma\le\frac{5(1+\nu)}{46+36\nu}-6\epsilon_0$  the following inequality holds
\begin{gather}\label{17-24}
\frac{1}{N}\mathop{{\sum}^*}_{1\le a\le q\le \xi_1^{1/(\nu+1)}\atop q>Q_0}
\int\limits_{q^{\nu}\le|K|\le\frac{\xi_1}{q}}
\left|S_N(\frac{a}{q}+\frac{K}{N})\right|^2dK\ll\frac{|\Omega_N|^2}{N}.
\end{gather}
\end{Le}
%%%%%%%%%%%%%%%%%%%%%%%%%%%%5
It remains to estimate the integral over the domain 5.
\begin{Le}\label{lemma-17-12}
For $\gamma\le\frac{1}{5+2\nu}-6\epsilon_0$  the following inequality holds
\begin{gather}\label{17-25}
\frac{1}{N}\mathop{{\sum}^*}_{1\le a\le q\le \xi_1\atop q>Q_0}
\int\limits_{\frac{Q_0}{q}\le|K|\le\min\{q^{\nu},\frac{\xi_1}{q}\}}
\left|S_N(\frac{a}{q}+\frac{K}{N})\right|^2dK
\ll\frac{|\Omega_N|^2}{N}.
\end{gather}
\begin{proof}
It follows immediately from the proof of Lemma 16.14 in ~\cite{FK} that for {$\gamma\le\frac{1}{6}-5\epsilon_0$} the following inequality holds
\begin{gather}\label{17-23}
\frac{1}{N}\mathop{{\sum}^*}_{1\le a\le q\le \xi_1\atop q>Q_0}
\int\limits_{\frac{Q_0}{q}\le|K|\le Q _0}
\left|S_N(\frac{a}{q}+\frac{K}{N})\right|^2dK\ll\frac{|\Omega_N|^2}{N}.
\end{gather}
Actually, we use Lemma \ref{lemma-17-5} with
\begin{gather*}
c_1=Q_0,\,c_2=\xi_1,\,
f_1=\frac{Q_0}{q},\,f_2=Q_0,\,
m_1=\frac{Q_0}{Q_1},\,m_2=Q_0.
\end{gather*}
It was proved in ~\cite[(16.45)]{FK} that
\begin{gather}\label{17-22-1}
\mathop{{\sum}^*}_{N_j\le q\le N_{j+1}\atop 1\le a\le q }\left|S_N(\frac{a}{q}+\frac{K}{N})\right|^2\ll
|\Omega_N|^2\frac{\KK^{4\gamma+12\epsilon_0}Q^{6\gamma+1+20\epsilon_0}}{\KK Q^2_1}.
\end{gather}
Integrating over $K$ and taking into account $m_2\ll1$ we obtain
\begin{gather}\label{17-23}
\int\limits_{m_1\le|K|\le m_2}
\mathop{{\sum}^*}_{N_j\le q\le N_{j+1}\atop 1\le a\le q }\left|S_N(\frac{a}{q}+\frac{K}{N})\right|^2\ll
|\Omega_N|^2Q^{6\gamma-1+22\epsilon_0}.
\end{gather}
For the sum over $j$ to be bounded by a constant, it is sufficient to have $\gamma\le\frac{1}{6}-5\epsilon_0.$
Hence, it remains to prove that
\begin{gather}\label{17-22-1}
\frac{1}{N}\mathop{{\sum}^*}_{1\le a\le q\le \xi_1\atop q>Q_0}
\int\limits_{Q_0\le|K|\le\min\{q^{\nu},\frac{\xi_1}{q}\}}
\left|S_N(\frac{a}{q}+\frac{K}{N})\right|^2dK\ll\frac{|\Omega_N|^2}{N}.
\end{gather}
Using \eqref{12-3-1} and arguments similar to Lemma \ref{lemma-17-5}, we obtain
\begin{gather}\notag
\frac{1}{N}\mathop{{\sum}^*}_{1\le a\le q\le \xi_1\atop q>Q_0}
\int\limits_{Q_0\le|K|\le\min\{q^{\nu},\frac{\xi_1}{q}\}}
\left|S_N(\frac{a}{q}+\frac{K}{N})\right|^2dK\ll
\frac{1}{TN}\mathop{{\sum}^*}_{1\le a\le q\le \xi_1\atop q>Q_0}
\sum_{TQ_0\le|l|\le T\min\{q^{\nu},\frac{\xi_1}{q}\}}
\left|S_N(\frac{a}{q}+\frac{l}{TN})\right|^2\le\\\ll
\frac{1}{TN}\sum_{j:\,Q_0^{1-\epsilon_0}\le N_j\le \xi_1\,}
\sum_{i:\,Q_0^{1-\epsilon_0}\le N_i\le \min\{N_{j+1}^{\nu},\frac{\xi_1}{N_j}\}\,}
\mathop{{\sum}^*}_{N_j\le q\le N_{j+1}\atop 1\le a\le q }
\sum_{TN_i\le|l|\le TN_{i+1}}
\left|S_N(\frac{a}{q}+\frac{l}{TN})\right|^2.\label{17-24}
\end{gather}
Applying Lemma \ref{lemma-15-1} with $Q_1=N_j, Q=N_{j+1}, \kappa_1=TN_i,\kappa=TN_{i+1}$ and taking into account that the number of summands in the sum over $i$ is less than $c\log\log N_{j+1}$ we obtain
\begin{gather}\notag
\frac{1}{N}\mathop{{\sum}^*}_{1\le a\le q\le \xi_1\atop q>Q_0}
\int\limits_{Q_0\le|K|\le\min\{q^{\nu},\frac{\xi_1}{q}\}}
\left|S_N(\frac{a}{q}+\frac{K}{N})\right|^2dK\ll\\\notag\ll
\frac{|\Omega_N|^2}{N}
\sum_{j:\,Q_0^{1-\epsilon_0}\le N_j\le \xi_1\,}
\sum_{i:\,Q_0^{1-\epsilon_0}\le N_i\le \min\{N_{j+1}^{\nu},\frac{\xi_1}{N_j}\}\,}
\left(N_{j+1}^2N_{i+1}\right)^{2\gamma+7\epsilon_0}N_{j}^{\gamma-1+5\epsilon_0}\ll\\\ll
\frac{|\Omega_N|^2}{N}
\sum_{j:\,Q_0^{1-\epsilon_0}\le N_j\le \xi_1\,}N_j^{(5+2\nu)\gamma-1+45\epsilon_0}.\label{17-25}
\end{gather}
For the sum over $j$ to be bounded by a constant, it is sufficient to have $\gamma\le\frac{1}{5+2\nu}-6\epsilon_0.$ This completes the proof of the lemma.
\end{proof}
\end{Le}
Having set $\nu=\frac{3}{2},$  we obtain that for $\gamma\le\frac{1}{8}-6\epsilon_0,\,\epsilon_0\in(0,\frac{1}{2500})$ the first integral in the right side of \eqref{17-4} is less than  $\frac{|\Omega_N|^2}{N}.$ So the inequality \eqref{result} holds for
$\gamma\le\frac{1}{8}-6\epsilon_0,\,\epsilon_0\in(0,\frac{1}{2500})$ and Theorem \ref{uslov2} is proved.\par
As mentioned in the paper ~\cite{FK2}, the proof of Lemma \ref{lemma-17-11} significantly uses the results of the paper ~\cite{BKS}. The following version of Lemma \ref{lemma-17-11} was proved in ~\cite{FK2} by elementary methods with the use of the estimates of Kloosterman sums.
%%%%%%%%%%%%%%%%%%%%%%%%%%%%%%%
\begin{Le}\label{lemma-17-11-2}(~\cite[лемма 8.10.]{FK2})
For $\gamma\le\frac{\nu-1/2}{10(1+\nu)}-8\epsilon_0$ the following inequality holds
\begin{gather}\label{17-24}
\frac{1}{N}\mathop{{\sum}^*}_{1\le a\le q\le \xi_1^{1/(\nu+1)}\atop q>Q_0}
\int\limits_{q^{\nu}\le|K|\le\frac{\xi_1}{q}}
\left|S_N(\frac{a}{q}+\frac{K}{N})\right|^2dK\ll\frac{|\Omega_N|^2}{N}.
\end{gather}
\end{Le}
%%%%%%%%%%%%%%%%%%%%%%%%%%%%5
Having set $\nu=\frac{3+\sqrt{34}}{2},$  we obtain that for $\gamma\le\frac{1}{8+\sqrt{34}}-8\epsilon_0,\,\epsilon_0\in(0,\frac{1}{2500})$ the first integral in the right side of \eqref{17-4} is less than $\frac{|\Omega_N|^2}{N}.$ So the inequality \eqref{result} holds for $\gamma\le\frac{1}{8+\sqrt{34}}-8\epsilon_0,\,\epsilon_0\in(0,\frac{1}{2500}).$ Hence, the following theorem is valid.
\begin{Th}\label{ucondition}
For any alphabet $\A$ with
\begin{equation}\label{KFcondition2}
\delta_{\A}>1-\frac{1}{8+\sqrt{34}}=0,9276
\end{equation}
the inequality \eqref{BKresult} holds.
\end{Th}
\begin{Zam}
It is proved ~\cite{Jenkinson} that $\delta_{10}=0,9257\ldots.$ From this follows that the alphabet $\left\{1,2,\ldots,10,11\right\}$ seems to satisfy the condition of Theorem \ref{ucondition}.
\end{Zam}

%%%%%%%%%%%%%%%%%%%%%%%%%%%%%%%%%%%%%%%%%%%%%%%%%%%%%%%%%%%%%%%%%%%%%%%%%%%%%%%%%%%%%%%%%%%%%%%%%%%%%%%%%%%%%%%%%%%%%%%%%%%%%%%%%%
%%%%%%%%%%%%%%%%%%%%%%%%%%%%%%%%%%%%%%%%%%%%%%%%%%%%%%%%%%%%%%%%%%%%%%%%%%%%%%%%%%%%%%%%%%%%%%%%%%%%%%%%%%%%%%%%%%%%%%%%%%%%%%%%%%

%\newpage

\end{document}